\documentclass[12pt]{article}
\usepackage{amsfonts}
\usepackage{mathrsfs}
\usepackage{amsmath,amssymb}
\baselineskip 4.5pt \lineskip 4.5pt \numberwithin{equation}{section}

\def\QED{\hfill$\Box$\par}

\def\HH{\mathcal {H}}

\def\l{\lambda}
\def\LL{\mathcal {L}_{\l,\mu}}

\def\cl{\centerline}

\def\ni{\noindent}

\def\rar{\longrightarrow}

\def\vs{\vspace*}

\def\C{\mathbb{C}}

\def\Z{\mathbb{Z}}
\def\adddot{$\!\!\!${\bf.}\ \ }

\newtheorem{theo}{Theorem}[section]

\newtheorem{case}{Case}
\newtheorem{lemm}[theo]{Lemma}

\begin{document}
\baselineskip 18pt
\cl{{\large \bf
 2-Cocycles of Deformative
Schr\"{o}dinger-Virasoro Algebras}\footnote {Supported by NSF grants
10471091, 10671027 of China, ``One Hundred Talents Program'' from
University of Science
and Technology of China.\\[2pt] \indent Corresponding E-mail:
sd\_junbo@163.com}} \vs{6pt}

\cl{Junbo Li$^{*,\dag)}$, Yucai Su$^{\ddag)}$}

\cl{\small $^{*)}$Department of Mathematics, Shanghai Jiao Tong
University, Shanghai 200240, China}

\cl{\small $^{\dag)}$Department of Mathematics, Changshu Institute
of Technology, Changshu 215500, China}

\cl{\small $^{\ddag)}$Department of Mathematics, University of
Science and Technology of China, Hefei 230026, China}

\cl{\small E-mail: sd\_junbo@163.com, ycsu@ustc.edu.cn}\vs{6pt}

\noindent{\small{\bf Abstract.} In a series of papers by Henkel,
Roger and Unterberger, Schr\"{o}dinger-Virasoro algebras and their
deformations were introduced and investigated. In the present paper
we determine the 2-cocycles of a class of deformative
Schr\"{o}dinger-Virasoro algebras.

} \vs{12pt}

\cl{\bf\S1. \
Introduction}\setcounter{section}{1}\setcounter{theo}{0}\setcounter{equation}{0}

It is well known that the infinite-dimensional Schr\"{o}dinger Lie
algebras and Virasoro algebra play important roles in many areas of
mathematics and physics (e.g., statistical physics). The
Schr\"{o}dinger-Virasoro algebras and their deformations were
introduced in \cite{H1,H2,HU,RU}, in the context of non-equilibrium
statistical physics, closely related to both Schr\"{o}dinger Lie
algebras and the Virasoro Lie algebra. Their vertex algebra
representations were constructed in \cite{U}, and later the
derivation algebra and the automorphism group of the twisted sector
were determined in \cite{LS1}. Furthermore, irreducible modules with
finite-dimensional weight spaces and indecomposable modules over
both original and twisted sectors were investigated in \cite{LS2}.

The infinite-dimensional Lie algebras $\LL\,(\l,\mu\in\C)$
considered in this paper called {\it twisted deformative
Schr\"{o}dinger-Virasoro Lie algebras} (\,see \cite{RU}), possess
the same $\C$-basis
$$\{L_n,\,M_n,\,Y_n\,|\,n\in \Z\}$$
with the following Lie brackets:
\begin{eqnarray}
&[L_n,L_m]\!\!\!&=(m-n)L_{m+n},\label{LB1}
\\[4pt]
&[L_n,Y_m]\!\!\!&=(m-\frac{(\l+1)n}{2}+\mu)Y_{m+n},\ \ \,\,
[Y_n,Y_{m}]=(m-n)M_{m+n},\label{LB2}\\[4pt]
&[L_n,M_m]\!\!\!&=(m-\l n+2\mu)M_{m+n},\ \ \ \,\ \ \ \ \,\,
[Y_n,M_m]=[M_n,M_m]=0.\label{LB3}
\end{eqnarray}

The purpose of this paper is to determine the 2-cocycles of
deformative Schr\"{o}dinger-Virasoro algebras $\LL\,(\l,\mu\in\C)$
defined above. The 2-cocycles on Lie algebras play important roles
in the central extensions of Lie algebras, which can be used to
construct many infinite-dimensional Lie algebras, such as affine Lie
algebras, Heisenberg algebras with a profound mathematical and
physical background, and further to describe the structures and some
of the representations of these Lie algebras. It is well known that
all 1-dimensional central extensions of some $\LL$ determine its
2-cohomology group. Since the cohomology groups are closely related
to the structure of Lie algebras, the computation of cohomology
groups seems to be important and interesting as well. Partially due
to the reasons stated above, there appeared a number of papers on
2-cocycles and cohomology groups of infinite-dimensional Lie
algebras and conformal algebras (\,see \cite{BKV}, \cite{L,LW},
\cite{SZR}--\cite{SZK} and related references cited in those
papers). Now let's formulate our main results below.

We start with a brief definition. Recall that a {\it 2-cocycle} on
some $\LL$ is a $\C$-bilinear function $\psi:\LL\times \LL\rar \C$
satisfying the following conditions:
\begin{eqnarray}
&&\psi(v_1,v_2)=-\psi(v_2,v_1)\mbox{\ \ (\,skew-symmetry)},\nonumber\\
&&\psi([v_1,v_2],v_3)+\psi([v_2,v_3],v_1)+\psi([v_3,v_1],v_2) =0
\mbox{\ \ (Jacobian identity)}\label{2c-Ji},
\end{eqnarray}
for $v_1,v_2,v_3\in \LL$. Denote by $\mathcal {C}^2(\LL,\C)$ the
vector space of 2-cocycles on $\LL$. For any $\C$-linear function
$f:\LL\rar\C$, one can define a 2-cocycle $\psi_f$ as follows
\begin{eqnarray}
\psi_f(v_1,v_2)=f([v_1,v_2]),\ \ \ \ \forall\;v_1,v_2\in
\LL.\label{cobco}
\end{eqnarray}
Such a 2-cocycle is called a {\it 2-coboundary} or a {\it trivial
2-cocycle} on $\LL$. Denote by $\mathcal {B}^2(\LL,\C)$ the vector
space of 2-coboundaries on $\LL$. A 2-cocycle $\varphi$ is said to
be {\it equivalent to} a 2-cocycle $\psi$ if $\varphi-\psi$ is
trivial. For a 2-cocycle $\psi$, we denote by $[\psi]$ the
equivalent class of $\psi$. The quotient space
\begin{eqnarray*}
\mathcal {H}^2(\LL,\C)\!=\!\mathcal {C}^2(\LL,\C)/\mathcal
{B}^2(\LL,\C)\! =\{\mbox{the equivalent classes of 2-cocycles}\},
\end{eqnarray*}
is called the {\it second cohomology group} of $\LL$.

For the case $\mu=0$, this problem has been considered and solved in
\cite{RU} by using the homological method. So we only need to
consider the case $\mu\in\C^*$ in the present paper. We call a
2-cocycle $\xi$ on $\LL$ the {\it Virasoro cocycle}, denoted by
$\xi_{Vir}$, if
\begin{eqnarray}\label{defvirco}
&&\xi(L_n,L_{m})=\frac{n^3-n}{12}\delta_{m,-n}\,,\ \ {\mbox {while
other components vanishing}}.
\end{eqnarray}
The main results of the paper can be formulated as the following
theorem.
\begin{theo}\adddot\label{mainth}
(i)\ If $\mu\notin\{\frac{1}{3}\Z\}$, then for any $\l\in\C$,
$\HH^2(\mathcal {L}{_{\l,\mu}},\C)\cong\C$ is generated
by the Virasoro cocycle.\\
(ii)\ If $\mu\in\{\frac{1}{3}+\Z,\,\frac{2}{3}+\Z\}$, then for any
$\l\ne-1$, also $\HH^2(\mathcal {L}{_{\l,\mu}},\C)\cong\C$ is
generated by the Virasoro cocycle.\\
(iii)\ If $\mu\in\{\frac{1}{3}+\Z,\,\frac{2}{3}+\Z\}$ and $\l=-1$,
$\HH^2(\mathcal {L}_{\l,\mu},\C)\cong\C^2$ is generated by the
Virasoro cocycle and an independent cocycle of the form
$c(M_{n},Y_m)=\delta_{n,-m-3\mu}$.\\
(iv)\ For $\mu\in\Z^*$ and $\l\ne-3,-1,1$, $\HH^2(\mathcal
{L}{_{\l,\mu}},\C)\cong\C$ is generated
by the Virasoro cocycle.\\
(v)\ For $\mu\in\Z^*$ and $\l=-3,1$, $\HH^2(\mathcal
{L}_{\l,\mu},\C)\cong\C^2$ is generated by the Virasoro cocycle and
an independent cocycle of the form
$c(L_n,Y_m)=\frac{m+\mu+1}{2}\delta_{n,-m-\mu}$ for $\l=-3$ or
$c(L_n,Y_m)=(m+\mu-1)(m+\mu)(m+\mu+1)\delta_{n,-m-\mu}$ for $\l=1$.\\
(vi)\ For $\mu\in\Z^*$ and $\l=-1$, $\HH^2(\mathcal
{L}_{\l,0},\C)\cong\C^3$ is generated by the Virasoro cocycle and
other two independent cocycles $c_1$ and $c_2$ defined by (all other
components vanishing)
\begin{eqnarray*}
c_1(L_n,Y_m)=\frac{(m+\mu)(m+\mu+1)}{2}\delta_{n,-m-\mu};\ \
c_2(M_n,Y_m)=\delta_{n,-m-3\mu}.
\end{eqnarray*}
\end{theo}

Throughout the paper, we denote by $\Z^*$ the set of all nonzero
integers, $\C^*$ the set of all nonzero complex numbers and
$\C^*\!\!\setminus\!\Z^*=\{x\,|\,x\in\C^*,\,x\notin\Z^*\}$.

\vs{18pt} \cl{\bf\S2. \ Proof the main
results}\setcounter{section}{2}\setcounter{theo}{0}\setcounter{equation}{0}

Let $\psi$ be any 2-cocycle. Our main object is to obtain all
equivalent classes of the nontrivial 2-cocycles by means of
subtracting all equivalent classes of the 2-coboundaries on $\LL$
from $\psi$.

Define a $\C$-linear function $f:\LL\rightarrow\C$ as follows
\begin{eqnarray}
f(L_n)\!\!\!&=&\!\!\!\left\{\begin{array}{ll}
\frac{1}{n}\psi(L_0,L_n) &{\rm if}\ \,n\neq0,\,\,\forall\,\,\mu\in\C^*,\vs{6pt}\\
\frac{1}{2}\psi(L_{-1},L_{1}) &{\rm if}\
\,n=0,\,\,\forall\,\,\mu\in\C^*,
\end{array}\right.\label{delfL}\\
f(M_n)\!\!\!&=&\!\!\!\left\{\begin{array}{ll}
\frac{1}{n+2\mu}\psi(L_0,M_n) &{\rm if}\
\,n\neq-2\mu,\,\mu\in\frac{1}{2}+\Z\cup\Z^*,\ {\rm or}\
\mu\notin\frac{1}{2}\Z,
\vs{6pt}\\
\frac{-1}{\l+1}\psi(L_1,M_{-2\mu-1}) &{\rm if}\
\,n=-2\mu,\,\l\neq-1, \ {\rm
and}\ \mu\in\frac{1}{2}+\Z\cup\Z^*,\end{array}\right.\label{delfM}\\
f(Y_n)\!\!\!&=&\!\!\!\left\{\begin{array}{ll}
\frac{1}{n+\mu}\psi(L_0,Y_n) &{\rm if}\ \,n\neq-\mu,\ \mu\in\Z^*,\
{\rm or}\ \mu\in\C^*\!\!\setminus\!\Z^*,\vs{6pt}\\
\frac{-2}{\l+3}\psi(L_1,Y_{-\mu-1}) &{\rm if}\ \,n=-\mu, \l\neq-3\
{\rm and}\ \mu\in\Z^*.\end{array}\right.\label{delfY}
\end{eqnarray}
Let $\varphi=\psi-\psi_f-\xi_{Vir}$ where $\psi_f$ and $\xi_{Vir}$
are respectively defined in (\ref{cobco}) and (\ref{defvirco}), then
\begin{eqnarray}\label{recl1}
\varphi(L_m,L_n)=0,\ \ \forall\,\,m,\,n\in\Z.
\end{eqnarray}
\begin{case}
\adddot\label{cs1} $\mu\notin\frac{1}{2}\Z$.
\end{case}
\begin{lemm}\adddot\label{lemm1}
If $\mu\in\{\frac{1}{3}+\Z,\,\frac{2}{3}+\Z\}\ \ {\rm and}\ \
\l=-1,$ one has (other components vanishing)
\begin{eqnarray}\label{lemconl1}
&&\varphi(M_{-m-3\mu},Y_m)=\varphi(M_{-3\mu},Y_0),\ \
\forall\,\,m,n\in\Z;
\end{eqnarray}
otherwise, $\varphi=0$.
\end{lemm}
{\it Proof.} According to (\ref{delfM}) and (\ref{delfY}), one has
\begin{eqnarray}\label{recl18998}
\varphi(L_0,Y_n)=\varphi(L_0,M_n)=0\ \ \forall\,\,n\in\Z.
\end{eqnarray}
For any $m,n\in\Z$, using the Jacobian identity on the triple
$(L_0,\,Y_{m},\,Y_n)$, together with (\ref{recl1}), we obtain
\begin{eqnarray*}
(m+n+2\mu)\varphi(Y_m,Y_n)=0,
\end{eqnarray*}
which together with our assumption $\mu\notin\frac{1}{2}\Z$ forces
\begin{eqnarray}\label{corI-4}
\varphi(Y_m,Y_n)=0.
\end{eqnarray}
Using the Jacobian identity on the four triples
$(L_{m},\,Y_n,\,L_0),\,(L_{m},\,M_n,\,L_0),\,(Y_{m},\,M_n,\,L_0)$
and $(M_{m},\,M_n,\,L_0)$ in (\ref{2c-Ji}) respectively, one has
\begin{eqnarray}
&&(\,m+n+\mu\,)\,\varphi(\,L_m,Y_n\,)\,=0,\label{corIad-1}\\
&&(m+n+2\mu)\,\varphi(L_m,M_n)=0,\label{corIad-2}\\
&&(m+n+3\mu)\,\varphi(Y_m,M_n)\,=0,\label{corI-5}\\
&&(m+n+4\mu)\varphi(M_m,M_n)=0.\label{corI-6}
\end{eqnarray}
Then our assumption $\mu\notin\frac{1}{2}\Z$, together with
(\ref{corIad-1}) and (\ref{corIad-2}) gives
\begin{eqnarray}
\varphi(L_m,Y_n)=\varphi(L_m,M_n)=0.\label{corIad-3}
\end{eqnarray}

\noindent{\bf Subcase 1.1.}\ \
$\mu\notin\{\frac{1}{3}+\Z,\,\frac{2}{3}+\Z,\,\frac{1}{4}+\Z,\,\frac{3}{4}+\Z\}$.

In this subcase, (\ref{corI-5}) and (\ref{corI-6}) force
\begin{eqnarray}
\varphi(Y_m,M_n)=\varphi(M_m,M_n)=0.\label{corI.1-1}
\end{eqnarray}

\noindent{\bf Subcase 1.2.}\ \
$\mu\in\{\frac{1}{3}+\Z\}\bigcup\{\frac{2}{3}+\Z\}$.

In this subcase, $m+n+4\mu\ne0$ for any $m,n\in\Z$. Then
(\ref{corI-5}) and (\ref{corI-6}) imply
\begin{eqnarray}
&&\varphi(\,Y_m,M_n\,)=0\ \ \ \ {\rm if}\
\,m+n\neq-3\mu,\label{corI.2-5-1}\\
&&\varphi(M_m,M_n)=0\ \ \ \ {\rm for\ any}\
\,m,\,n\in\Z.\label{corI.2-6-1}
\end{eqnarray}
We only need to compute the value of
\begin{eqnarray*}
\varphi(M_{-3\mu-m},Y_m)\ \ \,\,{\mbox{for any }}m\in\Z.
\end{eqnarray*}
Using the Jacobian identity on the triple
$(L_{-m},\,Y_{m},\,M_{-3\mu})$, one has
\begin{eqnarray}\label{corI.1ad-1}
&&\big((\l+3)m+2\mu\big)\varphi(M_{-3\mu},Y_0)+2(\l
m-\mu)\varphi(M_{-m-3\mu},Y_m)=0.
\end{eqnarray}
If $\l=0$, then (\ref{corI.1ad-1}) gives (\,since $\mu\neq0$ in this
case)
\begin{eqnarray}\label{corI.1ad-1.1}
&&\varphi(M_{-m-3\mu},Y_m)=(\frac{3m}{2\mu}+1)\varphi(M_{-3\mu},Y_0).
\end{eqnarray}
Applying the Jacobian identity on the triple
$(L_{-m-3\mu},Y_{m},M_{0})$ for any $m\in\Z$, together with
(\ref{corI.1ad-1.1}), one has $(m+3\mu)\varphi(M_{-3\mu},Y_0)=0$,
which gives $\varphi(M_{-3\mu},Y_0)=0$, and further \big(together
with (\ref{corI.1ad-1.1})\big)
\begin{eqnarray}\label{yzI.1ad-1.2}
\varphi(M_{-m-3\mu},Y_m)=0.
\end{eqnarray}
If $\l\in\C^*$ and $\frac{\mu}{\l}\notin\Z$, then (\ref{corI.1ad-1})
gives
\begin{eqnarray}\label{corI.1ad-1.2}
&&\varphi(M_{-m-3\mu},Y_m)=\frac{(\l+3)m+2\mu}{2(\mu-m\l)}\varphi(M_{-3\mu},Y_0).
\end{eqnarray}
Applying the Jacobian identity on the triple
$(L_{-m-3\mu},Y_{m},M_{0})\,(\,\forall\,\,m\in\Z)$, together with
(\ref{corI.1ad-1.2}), one has (\,our assumption forcing
$\l\ne-\frac{1}{3}$)
\begin{eqnarray*}
\frac{9(1+\l)(m+3\mu)\big(m
\l(3+\l)+\mu(\l-1)\big)}{4(1+3\l)(m\l-mu)}\varphi(M_{-3\mu},Y_0)=0,
\end{eqnarray*}
which gives $\varphi(M_{-3\mu},Y_0)=0$ if $\l\ne-1$ and further
\big(together with (\ref{corI.1ad-1.2})\big)
\begin{eqnarray}\label{yzI.1ad-1.2}
\varphi(M_{-m-3\mu},Y_m)=0\ \ \ \,{\mbox {if}\ }\l\ne-1.
\end{eqnarray}
For the case $\l=-1$, (\ref{corI.1ad-1.2}) gives
\begin{eqnarray}\label{lastI.1ad-1.2}
&&\varphi(M_{-m-3\mu},Y_m)=\varphi(M_{-3\mu},Y_0).
\end{eqnarray}
If $\l\in\C^*$ while $\frac{\mu}{\l}\in\Z$, then (\ref{corI.1ad-1})
gives (\,since$\l\ne-1,-3$ in this subcase)
\begin{eqnarray}
&&(1+\frac{1}{\l})\mu\varphi(M_{-3\mu},Y_0)=0\ \ {\rm and}\ \
\varphi(M_{-3\mu},Y_0)=0,\nonumber\\
&&\varphi(M_{-m-3\mu},Y_m)=\frac{(\l+3)m+2\mu}{2(\mu-m\l)}\varphi(M_{-3\mu},Y_0)=0\
\ \ {\rm if}\ \ m\ne\frac{\mu}{\l}.\label{corI.1ad-1.4}
\end{eqnarray}
Using the Jacobian identity on
$(L_{1},\,Y_{\frac{\mu}{\l}-1},\,M_{-\frac{\mu}{\l}-3\mu})$ and
$(L_{2},\,Y_{\frac{\mu}{\l}-2},\,M_{-\frac{\mu}{\l}-3\mu})$,
together with (\ref{corI.1ad-1.4}), one has
\begin{eqnarray*}
&&(\l+3)\varphi(M_{-\frac{\mu}{\l}-3\mu},Y_{\frac{\mu}{\l}})=0\ \
{\rm i.e.,}\ \
\varphi(M_{-\frac{\mu}{\l}-3\mu},Y_{\frac{\mu}{\l}})=0\,\,({\rm as}\
\,\l\ne-3),
\end{eqnarray*}
which combined with (\ref{corI.1ad-1.4}) gives
\begin{eqnarray}
\varphi(M_{-3\mu-m},Y_m)=0,\ \ \,\,\forall\,\,m\in\Z.
\end{eqnarray}

\noindent{\bf Subcase 1.3.}\ \
$\mu\in\{\frac{1}{4}+\Z\}\bigcup\{\frac{3}{4}+\Z\}$.

In this subcase, $m+n+3\mu\ne0$ for any $m,n\in\Z$. Then
(\ref{corI-5}) and (\ref{corI-6}) give
\begin{eqnarray}
&&\varphi(M_m,M_n)=0\ \ \ \ {\rm if}\
\,m+n\neq-4\mu,\label{corcoI.2-5-1}\\
&&\varphi(\,Y_m,M_n\,)=0\ \ \ \ {\rm for\ any}\
\,m,\,n\in\Z.\label{corI.2-5-1}
\end{eqnarray}
Then the left sector we have to consider in this subcase is
$\varphi(M_{-4\mu-m},M_m)$ for any $m\in\Z$.

Using the Jacobian identity on the triple
$(Y_{-4\mu},\,Y_{-m},\,M_m)$, one has
\begin{eqnarray*}
&&-(m-4\mu)\varphi(M_{-m-4\mu},M_m)=0,
\end{eqnarray*}
which gives
\begin{eqnarray}
&&\varphi(M_{-m-4\mu},M_m)=0\ \ \,{\rm for\ \,any}\ \ 4\mu\ne
m\in\Z.\label{I-IMMa1-1}
\end{eqnarray}
Using the Jacobian identity on the triples
$(Y_{-6\mu},\,Y_{-2\mu},\,M_{4\mu})$, one has
\begin{eqnarray*}
4\mu\varphi(M_{-8\mu},M_{4\mu})=0,
\end{eqnarray*}
which gives $\varphi(M_{-8\mu},M_{4\mu})=0$ and further together
with (\ref{corcoI.2-5-1}) gives
\begin{eqnarray}
&&\varphi(M_m,M_n)=0,\ \ \forall\,\,m,n\in\Z.\label{corI-I-3}
\end{eqnarray}
Then this lemma follows.\QED

\ni This lemma in particular proves Theorem \ref{mainth} (i)--(iii)
in this case.
\begin{case}
\adddot\label{cs2} $\mu\in\frac{1}{2}+\Z$.
\end{case}
\begin{lemm}\adddot\label{lemm2}
In this case $\varphi(x,y)=0,\,\ \forall\,\,\l\in\C,\,x,y\in\LL$.
\end{lemm}
{\it Proof.} According to (\ref{delfM}) and (\ref{delfY}), one has
\begin{eqnarray}
&&\varphi(L_0,Y_n)=0\ \,\,\forall\,\,n\in\Z,\label{recII2}\\
&&\varphi(L_0,M_n)=0\ \ {\rm if}\ \,n\neq-2\mu,\ \
\varphi(L_1,M_{-2\mu-1})=0\ \ {\rm if}\ \,\l\neq-1.\label{recII1}
\end{eqnarray}
Using the Jacobian identity by replacing the triple $(v_1,v_2,v_3)$
by five triples $(L_{m},\,Y_n,\,L_0)$, $(Y_{m},\,M_n,\,L_0)$,
$(M_{m},\,M_n,\,L_0)$, $(Y_{m},\,Y_n,\,L_0)$ and
$(L_{m},\,M_n,\,L_0)$ in (\ref{2c-Ji}) respectively, together with
(\ref{recII2}), one has
\begin{eqnarray}
&&(m+n+\mu)\varphi(L_m,Y_n)=0,\label{corIILY}\\
&&(m+n+3\mu)\varphi(Y_m,M_n)=0,\label{corIIYM}\\
&&(m+n+4\mu)\varphi(M_m,M_n)=0,\label{corIIMM}\\
&&(m-n)\varphi(L_0,M_{m+n})+(m+n+2\mu)\varphi(Y_m,Y_n)=0,\label{corIIYY}\\
&&(-n+m\l-2\mu)\varphi(L_0,M_{m+n})+(m+n+2\mu)\varphi(L_m,M_n)=0.\label{corIILM}
\end{eqnarray}
In this case, $(m+n+\mu)(m+n+3\mu)\ne0$. Hence (\ref{corIILY}) and
(\ref{corIIYM}) force
\begin{eqnarray}
\varphi(L_m,Y_n)=\varphi(Y_m,M_n)=0.\label{corIIYMLY}
\end{eqnarray}
According to (\ref{recII1}), (\ref{corIIYY}) and (\ref{corIILM}),
one can deduce
\begin{eqnarray}
&&\varphi(\,Y_m,Y_n\,)=0\ \ \ {\rm if}\
\,m+n\neq-2\mu,\label{corIIYYl2}\\
&&\varphi(L_m,M_n)=0\ \ \ {\rm if}\
\,m+n\neq-2\mu,\label{corIIYYl30}\\
&&-2(n+\mu)\varphi(L_0,M_{-2\mu})=0\ \ \ {\rm if}\
\,m+n=-2\mu.\label{corIIYYl1}
\end{eqnarray}
The assumption $\mu\in\frac{1}{2}+\Z$ and (\ref{corIIYYl1}) infer
\begin{eqnarray}
\varphi(L_0,M_{-2\mu})=0.\label{corIIYYl1a}
\end{eqnarray}
From (\ref{recII1}), (\ref{recII2}), (\ref{corIIMM}),
(\ref{corIIYMLY}), (\ref{corIIYYl2}) and (\ref{corIIYYl1a}), the
left components we have to present in this case are listed in the
following (where $m$ is an arbitrary integer):
\begin{eqnarray*}
\varphi(Y_{-2\mu-m},Y_m),\ \ \varphi(L_{-2\mu-m},M_m)\ \ {\rm and}\
\ \varphi(M_{-4\mu-m},M_m).
\end{eqnarray*}
Repeating the proving process between (\ref{corcoI.2-5-1}) and
(\ref{corI-I-3}), in this case one also can obtain
\begin{eqnarray}
&&\varphi(M_{-4\mu-m},M_m)=0,\ \ \forall\,\,m\in\Z.\label{corII-I-1}
\end{eqnarray}

Using the Jacobian identity on $(L_{-m},\,Y_{-2\mu},\,Y_m)$ and
$(L_{-2\mu-m},\,L_{n},\,M_{m-n})\,\,(\,\forall\,\,n\in\Z)$
respectively, one has
\begin{eqnarray}
&&\!\!\!\!\!\!\!\!\!\big(m(1+\l)-2\mu\big)\varphi(Y_{-m-2\mu},Y_{m})\nonumber\\
&&\!\!\!\!\!\!\!\!\!=2(m+2\mu)\varphi(L_{-m},M_{m-2\mu})-\big(m(3+\l)+2\mu\big)
\varphi(Y_{-2\mu},Y_0),\label{corII.2lad-1}\\
&&\!\!\!\!\!\!\!\!\!\big(m+2\mu-n(1+\l)\big)\varphi(L_{-m-2\mu},M_{m})\nonumber\\
&&\!\!\!\!\!\!\!\!\!=\big((m+2\mu)(1+\l)-n\big)\varphi(L_{n},M_{-n-2\mu})+(m+2\mu+n)
\varphi(L_{n-m-2\mu},M_{m-n}).\label{corII.1ad-n}
\end{eqnarray}
Replacing $n$ by $-n$ and $m$ by $m+n$ in (\ref{corII.1ad-n}), one
respectively gets
\begin{eqnarray}
&&\!\!\!\!\!\!\!\!\!\!\!\!\!\!\!\!\!\!\!
\big(m+2\mu+n(1+\l)\big)\varphi(L_{-m-2\mu},M_{m})\nonumber\\
&&\!\!\!\!\!\!\!\!\!\!\!\!\!\!\!\!\!\!\!
=(m+2\mu-n)\varphi(L_{-n-m-2\mu},M_{m+n})+
\big((m+2\mu)(1+\l)+n\big)\varphi(L_{-n},M_{n-2\mu}),\label{corII.1adn}\\
&&\!\!\!\!\!\!\!\!\!\!\!\!\!\!\!\!\!\!\!
(2n+m+2\mu) \varphi(L_{-m-2\mu},M_{m})\nonumber\\
&&\!\!\!\!\!\!\!\!\!\!\!\!\!\!\!\!\!\!\!
=(m+2\mu-n\l)\varphi(L_{-n-m-2\mu},M_{m+n})
-\big((n+m+2\mu)(1+\l)-n\big)\varphi(L_{n},M_{-n-2\mu}).\label{corII.1adm+n}
\end{eqnarray}

According to (\ref{recII1}), we have to divide the left part of this
proof into two subcases.

\vs{8pt}

\noindent{\bf Subcase 2.1.}\ \ $\l\ne-1$.

\vs{6pt}

In this subcase, taking $n=-1$ in both (\ref{corII.1adn}) and
(\ref{corII.1adm+n}), together with (\ref{recII1}), we have
\begin{eqnarray}
&&\!\!\!\!\!\!\!\!\!\!\!\!\!(m+2\mu-1-\l)\varphi(L_{-m-2\mu},M_{m})
=(m+2\mu+1)\varphi(L_{1-m-2\mu},M_{m-1}),\label{corII.1ad-2-3}\\
&&\!\!\!\!\!\!\!\!\!\!\!\!\!(m+2\mu-2)\varphi(L_{-m-2\mu},M_{m})\nonumber\\
&&\!\!\!\!\!\!\!\!\!\!\!\!\!=(m+2\mu+\l)
\varphi(L_{1-m-2\mu},M_{m-1})-\big((m+2\mu-1)(1+\l)+1\big)
\varphi(L_{-1},M_{1-2\mu}).\label{corII.1adm-1}
\end{eqnarray}
If $({\l}^2+\l-2)(1+\l)\ne0$, then combining (\ref{corII.1ad-2-3})
with (\ref{corII.1adm-1}), one can deduce
\begin{eqnarray}\label{reII1}
\varphi(L_{-m-2\mu},M_m)
=-\frac{(m+2\mu+1)\big((m+2\mu)(1+\l)-\l\big)}{{\l}^2+\l-2}
\varphi(L_{-1},M_{1-2\mu}).
\end{eqnarray}
Taking $\varphi(L_{-m-2\mu},M_m)$ and
$\varphi(L_{1-m-2\mu},M_{m-1})$ obtained from (\ref{reII1}) back to
(\ref{corII.1ad-2-3}), one has
\begin{eqnarray*}
\frac{\l(1+\l)\big((m+2\mu)^2-1\big)\varphi(L_{-1},M_{1-2\mu})}{{\l}^2+\l-2}=0,
\end{eqnarray*}
which forces (\,since the index $m$ can be shifted and our
assumption $\l\ne-1$)
\begin{eqnarray*}
\l\varphi(L_{-1},M_{1-2\mu})=0.
\end{eqnarray*}
In another word,
\begin{eqnarray}\label{cali1}
&&\!\!\!\!\!\!\!\!\!\!\!\!{\mbox{the system consisted of linear
equations
(\ref{corII.1ad-2-3}) and (\ref{corII.1adm-1}) has nonzero}}\nonumber\\
&&\!\!\!\!\!\!\!\!\!\!\!\!{\mbox{solutions if and only if $\l=0$
under our assumption $({\l}^2+\l-2)(1+\l)\ne0$.}}
\end{eqnarray}

Hence based on the discussions between (\ref{corII.1ad-2-3}) and
(\ref{cali1}), we have to divide Subcase 2.1 into another four
subcases.

\vs{8pt}

\noindent{\bf Subcase 2.1(i).}\ \ $\l=0$.

\vs{6pt}

If $\l=0$, then (\ref{reII1}) can be rewritten as
\begin{eqnarray}\label{reII1l=066}
\varphi(L_{-m-2\mu},M_m) =\frac{(m+2\mu)(m+2\mu+1)}{2}
\varphi(L_{-1},M_{1-2\mu}).
\end{eqnarray}

If $\l=0$ and $m\ne2\mu$, then (\ref{reII1l=066}) together with
(\ref{corII.2lad-1}) gives
\begin{eqnarray}\label{corII.2lad-3}
\varphi(Y_{-m-2\mu},Y_{m})=\frac{m(m+1)(m+2\mu)}{m-2\mu}
\varphi(L_{-1},M_{1-2\mu})-\frac{3m+2\mu}{m-2\mu}\varphi(Y_{-2\mu},Y_0).
\end{eqnarray}
Using the Jacobian identity on the triple
$(L_{-2\mu+m},\,Y_{0},\,Y_{-m})$, one has
\begin{eqnarray}\label{corII.2nesp2}
&&\big(2u(2+\l)-m(1+\l)\big)\varphi(Y_{m-2\mu},Y_{-m})\nonumber\\
&&=\big(m(3+\l)-2\mu(2+\l)\big)\varphi(Y_{0},Y_{-2\mu})-2m\varphi(L_{m-2\mu},M_{-m}).
\end{eqnarray}
Taking $\l=0,\,m=-2\mu$ in (\ref{corII.2nesp2}) and using
(\ref{reII1l=066}) together with (\ref{corII.2lad-3}), one has
\begin{eqnarray}\label{corII.2nesp2-1}
&&6\mu\varphi(Y_{-4\mu},Y_{2\mu})
=10\mu\varphi(Y_{-2\mu},Y_{0})+8\mu^2(4\mu+1)\varphi(L_{-1},M_{1-2\mu}),
\end{eqnarray}
which gives
\begin{eqnarray}\label{corII.2nesp2-166}
\varphi(Y_{-4\mu},Y_{2\mu})
=\frac{5}{3}\varphi(Y_{-2\mu},Y_{0})+\frac{4\mu(4\mu+1)}{3}\varphi(L_{-1},M_{1-2\mu}).
\end{eqnarray}
For any $p\in\Z^*$, applying the Jacobian identity on the triple
$(L_{p},Y_{2\mu},Y_{-p-4\mu})$, together with (\ref{reII1l=066}),
(\ref{corII.2lad-3}) and (\ref{corII.2nesp2-166}), we obtain
\begin{eqnarray*}
(p+4\mu)\Big((p+6\mu)\varphi(Y_{-2\mu},Y_0)
+\big(p^3-6(1+2\mu){\mu}^2-p\mu(1+6\mu)\big)\varphi(L_{-1},M_{1-2\mu})\Big)=0,
\end{eqnarray*}
which forces
\begin{eqnarray*}
\varphi(Y_{-2\mu},Y_0)=\varphi(L_{-1},M_{1-2\mu})=0.
\end{eqnarray*}
and further \big(recalling (\ref{corII.2lad-3}) and
(\ref{corII.2nesp2-166})\big)
\begin{eqnarray}\label{cheII1}
\varphi(Y_{-m-2\mu},Y_{m})=0,\ \ \ \forall\,\,m\in\Z.
\end{eqnarray}

\vs{8pt}

\noindent{\bf Subcase 2.1(ii).}\ \ $\l=1$.

\vs{6pt}

If $\l=1$, then (\ref{corII.1ad-2-3}) becomes
\begin{eqnarray}\label{corl=1}
&&\!\!\!\!\!\!\!\!\!\!\!\!\!(m+2\mu-2)\varphi(L_{-m-2\mu},M_{m})
=(m+2\mu+1)\varphi(L_{1-m-2\mu},M_{m-1}),
\end{eqnarray}
which further gives \big(\,by taking $m=2-2\mu$ in
(\ref{corl=1})\big)
\begin{eqnarray}\label{reII2}
\varphi(L_{-1},M_{1-2\mu})=0.
\end{eqnarray}
Also by (\ref{recII1}), (\ref{corIIYYl1a}), (\ref{corl=1}) and
(\ref{reII2}), one can deduce
\begin{eqnarray}\label{reII.1-00a}
\varphi(L_{-m-2\mu},M_{m})=\left\{\begin{array}{ll}
(m+2\mu-1)(m+2\mu)(m+2\mu+1)c_1&{\rm if}\ \,m\geq-2\mu-1,\vs{6pt}\\
(m+2\mu-1)(m+2\mu)(m+2\mu+1)c_2 &{\rm if}\ \,m<-2\mu-1,
\end{array}\right.
\end{eqnarray}
for some constants $c_1,\,c_2\in\C$. One thing left to be done is to
find the relations between the constants $c_1$ and $c_2$. If $\l=1$,
then (\ref{corII.1adn}) and (\ref{corII.1adm+n}) become
\begin{eqnarray*}
&&\!\!\!\!\!\!\!\!\!\!\!\!\!(m+2\mu+2n)\varphi(L_{-m-2\mu},M_{m})\nonumber\\
&&\!\!\!\!\!\!\!\!\!\!\!\!\!=(m+2\mu-n)
\varphi(L_{-n-m-2\mu},M_{m+n})+\big(2(m+2\mu)+n\big)
\varphi(L_{-n},M_{n-2\mu}),\\
&&\!\!\!\!\!\!\!\!\!\!\!\!\!(m+2\mu+2n)\varphi(L_{-m-2\mu},M_{m})\nonumber\\
&&\!\!\!\!\!\!\!\!\!\!\!\!\!=(m+2\mu-n)\varphi(L_{-n-m-2\mu},M_{m+n})
-\big(2(m+2\mu)+n\big)\varphi(L_{n},M_{-n-2\mu}),
\end{eqnarray*}
which together with each other force
\begin{eqnarray}
\big(2(m+2\mu)+n\big)
\big(\varphi(L_{-n},M_{n-2\mu})+\varphi(L_{n},M_{-n-2\mu})\big)=0,
\end{eqnarray}
and in particular give (\,by taking $n=2$)
\begin{eqnarray}
&&\big(m+2\mu+1\big)
\big(\varphi(L_{-2},M_{2-2\mu})+\varphi(L_{2},M_{-2-2\mu})\big)=0.\label{reII.1-1a}
\end{eqnarray}
Noticing $2-2\mu>-1-2\mu,\,-2-2\mu<-1-2\mu$, and then combining
(\ref{reII.1-00a}) with (\ref{reII.1-1a}), one can safely deduce
\begin{eqnarray*}
c_1=c_2,
\end{eqnarray*}
which together with (\ref{reII.1-00a}) gives
\begin{eqnarray}\label{reII.1-las}
\varphi(L_{-m-2\mu},M_{m})=(m+2\mu-1)(m+2\mu)(m+2\mu+1)c_1,\ \
\forall\,\,m\in\Z.
\end{eqnarray}

If $\l=1$, then $m\ne\mu\,\,(\,\forall\,\,m\in\Z)$ and
(\ref{reII.1-las}) together with (\ref{corII.2lad-1}) gives
\begin{eqnarray}\label{corII.2lad-5}
&&\varphi(Y_{-m-2\mu},Y_{m})=\frac{\mu+2m}{\mu-m}\varphi(Y_{-2\mu},Y_0).
\end{eqnarray}
For any $p\in\Z^*$, applying the Jacobian identity on the triple
$(L_p,Y_{2\mu},Y_{-p-4\mu})$, together with (\ref{reII.1-las}) and
(\ref{corII.2lad-5}), one has
\begin{eqnarray*}
p(p+4\mu)\varphi(Y_{-2\mu},Y_0)-5c_1(p^2-1)(p^2+11p\mu+30{\mu}^2)=0,
\end{eqnarray*}
which implies $\varphi(Y_{-2\mu},Y_0)=c_1=0$ and further
\begin{eqnarray}\label{lachI1}
\varphi(L_{-m-2\mu},M_{m})=\varphi(Y_{-m-2\mu},Y_{m})=0.
\end{eqnarray}

\noindent{\bf Subcase 2.1(iii).}\ \ $\l=-2$.

\vs{6pt}

If $\l=-2$, then (\ref{corII.1adn}) and (\ref{corII.1adm+n}) convert
to the following form:
\begin{eqnarray}
&&\!\!\!\!\!\!\!\!\!\!\!\!(m+2\mu-n)\big(\varphi(L_{-m-2\mu},M_{m})
-\varphi(L_{-n-m-2\mu},M_{m+n})\big)\nonumber\\
&&\!\!\!\!\!\!\!\!\!\!\!\!=
(n-m-2\mu)\varphi(L_{-n},M_{n-2\mu}),\label{corIIl--2.1}\\
&&\!\!\!\!\!\!\!\!\!\!\!\!(m+2\mu+2n)\big(\varphi(L_{-m-2\mu},M_{m})
-\varphi(L_{n-m-2\mu},M_{m+n})\big)\nonumber\\
&&\!\!\!\!\!\!\!\!\!\!\!\!=(m+2\mu+2n)\varphi(L_{n},M_{-n-2\mu}).\label{corIIl--2.2}
\end{eqnarray}
Furthermore, taking $n=-1$ in both (\ref{corIIl--2.1}) and
(\ref{corIIl--2.2}), and using (\ref{recII1}), one has
\begin{eqnarray*}
(m+2\mu+1)\big(\varphi(L_{-m-2\mu},M_{m})-
\varphi(L_{1-m-2\mu},M_{m-1})\big)\!\!\!&=\!\!\!&0,\\
(m+2\mu-2)\big(\varphi(L_{-m-2\mu},M_{m})
-\varphi(L_{1-m-2\mu},M_{m-1})\big)\!\!\!&=\!\!\!&
(m+2\mu-1)\varphi(L_{-1},M_{1-2\mu}),
\end{eqnarray*}
from which and using (\ref{recII1}) again, can we deduce the
following relation:
\begin{eqnarray}\label{relalast0}
\varphi(L_{-m-2\mu},M_{m})=0,\ \ \forall\,\,m\in\Z.
\end{eqnarray}

If $\l=-2$ and $m\ne-2\mu$, then (\ref{relalast0}) together with
(\ref{corII.2lad-1}) gives
\begin{eqnarray}\label{corII.2lad-4}
&&\varphi(Y_{-m-2\mu},Y_{m})=\varphi(Y_{-2\mu},Y_0).
\end{eqnarray}
Using the Jacobian identity on the triple
$(L_{-2\mu},\,Y_{m},\,Y_{-m})$, one has
\begin{eqnarray}\label{corII.2nesp1}
&&\big(\mu(2+\l)+m\big)\varphi(Y_{m-2\mu},Y_{-m})\nonumber\\
&&=-2m\varphi(L_{-2\mu},M_0)-\big(\mu(2+\l)-m\big)
\varphi(Y_{m},Y_{-m-2\mu}).
\end{eqnarray}
Taking $\l=-2,\,m=2\mu$ in (\ref{corII.2nesp1}) and using
(\ref{relalast0}) together with (\ref{corII.2lad-4}), one has
\begin{eqnarray}\label{corII.2nesp12}
&&\varphi(Y_{0},Y_{-2\mu})=\varphi(Y_{2\mu},Y_{-4\mu})
=\varphi(Y_{-2\mu},Y_0),
\end{eqnarray}
which gives
\begin{eqnarray}\label{corIcorrc1}
&&\varphi(Y_{-m-2\mu},Y_{m})=\varphi(Y_{-2\mu},Y_0)=0,\ \ \
\forall\,\, m\in\Z.
\end{eqnarray}

\vs{8pt}

\noindent{\bf Subcase 2.1(iv).}\ \ $\l\notin\{-2,\,-1,\,0,\,1\}$.

\vs{6pt}

If $\l\notin\{-2,\,-1,\,0,\,1\}$, then
$\varphi(L_{-m},M_{m-2\mu})=0\,\,(\,\forall\,\,m\in\Z)$ and
(\ref{corII.2lad-1}) can be rewritten as
\begin{eqnarray}\label{corII.2lad-6}
&&\big(m(\l+1)-2\mu\big)\varphi(Y_{-m-2\mu},Y_{m})
=-\big(m(\l+3)+2\mu\big)\varphi(Y_{-2\mu},Y_0).
\end{eqnarray}
If $\frac{2\mu}{\l+1}\notin\Z$, then
\begin{eqnarray}\label{corII.2lad-6-1}
&&\varphi(Y_{-m-2\mu},Y_{m})
=-\frac{m(\l+3)+2\mu}{m(\l+1)-2\mu}\varphi(Y_{-2\mu},Y_0).
\end{eqnarray}
For $m=\frac{2\mu}{\l+1}\in\Z$, (\ref{corII.2lad-6}) gives (since
$\l\notin\{-2,\,-1,\,0,\,1\}$ in this case)
\begin{eqnarray}\label{corII.2lad-7}
&&\varphi(Y_{-2\mu},Y_0)=0.
\end{eqnarray}
For $m\ne\frac{2\mu}{\l+1}\in\Z$, then
\begin{eqnarray}\label{corII.2lad-6-1}
&&\varphi(Y_{-m-2\mu},Y_{m})
=\frac{m(\l+3)+2\mu}{\mu\big(m(\l+1)-2\mu\big)}\varphi(Y_{-2\mu},Y_0)=0.
\end{eqnarray}

For $\l\notin\{-2,\,-1,\,0,\,1\}$, taking $m=-\frac{2\mu}{\l+1}$ in
(\ref{corII.2nesp1}) and using (\ref{cali1}) together with
(\ref{corII.2lad-6-1}), one has
\begin{eqnarray*}
&&\frac{\mu(\l+1)(\l+2)-2\mu}{\l+1}
\varphi(Y_{-\frac{2\mu}{\l+1}-2\mu},Y_{\frac{2\mu}{\l+1}})\nonumber\\
&&=\frac{4\mu}{\l+1}\varphi(L_{-2\mu},M_0)-\frac{\mu(\l+1)(\l+2)+2\mu}{\l+1}
\varphi(Y_{-\frac{2\mu}{\l+1}},Y_{\frac{2\mu}{\l+1}-2\mu})\nonumber\\
&&=0,
\end{eqnarray*}
which gives (\,since $\l\ne-3$ under our assumption
$\frac{2\mu}{\l+1}\in\Z$)
\begin{eqnarray}\label{corII.2nesp3ad}
\varphi(Y_{-\frac{2\mu}{\l+1}-2\mu},Y_{\frac{2\mu}{\l+1}})=0.
\end{eqnarray}

\vs{8pt}

\noindent{\bf Subcase 2.2.}\ \ $\l=-1$.

\vs{6pt}

If $\l=-1,\,n=1$, then (\ref{corII.1ad-n}) can be rewritten as
\begin{eqnarray*}
(m+2\mu)\varphi(L_{-m-2\mu},M_{m})=(1+m+2\mu)\varphi(L_{1-m-2\mu},M_{m-1})-
\varphi(L_{1},M_{-1-2\mu}),
\end{eqnarray*}
which is equivalent to
\begin{eqnarray}\label{corII.1ad-2-4}
\varphi(L_{-m-2\mu},M_{m})=\left\{\begin{array}{ll}
(m+2\mu)\psi(L_1,M_{-1-2\mu}) &{\rm if}\ \,m\geq-2\mu,\vs{8pt}\\
(m+2\mu+2)\varphi(L_1,M_{-1-2\mu})\\
-(m+2\mu+1)\varphi(L_2,M_{-2-2\mu})&{\rm if}\ \,m<-2\mu.
\end{array}\right.
\end{eqnarray}

If $\l=-1$, then (\ref{corII.2lad-1}) becomes
\begin{eqnarray*}
&&\varphi(Y_{-m-2\mu},Y_{m})=\big(\frac{m}{\mu}+1\big)
\varphi(Y_{-2\mu},Y_0) -(\frac{m}{\mu}+2)\varphi(L_{-m},M_{m-2\mu}),
\end{eqnarray*}
which combined with (\ref{corII.1ad-2-4}), can be simplified as
\begin{eqnarray}\label{corII.2ad-2}
&&\varphi(L_{-m-2\mu},Y_{m})\nonumber\\
&&=\left\{\begin{array}{ll}
\big(\frac{m}{\mu}+1\big)\varphi(Y_{-2\mu},Y_0)\\
-m(\frac{m}{\mu}+2)\psi(L_1,M_{-1-2\mu}) &{\rm if}\ \,m\geq0,\vs{8pt}\\
\big(\frac{m}{\mu}+1\big) \varphi(Y_{-2\mu},Y_0)
-(\frac{m}{\mu}+2)\big((m+2)\\
\times\varphi(L_1,M_{-1-2\mu})-(m+1)\varphi(L_2,M_{-2-2\mu})\big)&{\rm
if}\ \,m<0.
\end{array}\right.
\end{eqnarray}
However, (\ref{corII.1ad-2-4}) and (\ref{corII.2ad-2}) are not
compatible with the Jacobian identity given in (\ref{2c-Ji}), which
forces both of them must be zero. \ni This completes the proof of
the lemma.\QED

\ni Then the lemma proves Theorem \ref{mainth} (i) in this case.
\begin{case}
\adddot\label{cs4} $\mu\in\Z^*$.
\end{case}
\begin{lemm}\adddot\label{lemm4}
(i)\ For the subcase $\l\ne-3,-1,1$, one has $\varphi=0$.\\
(ii)\ For the subcase $\l=-3$, only
$\varphi(L_n,Y_{-n-\mu})\,\,(\,\forall\,\,n\in\Z)$ is not
vanishing, given in (\ref{reIIILY1-2}).\\
(iii)\ For the subcase $\l=1$, only
$\varphi(L_n,Y_{-n-\mu})\,\,(\,\forall\,\,n\in\Z)$ is not vanishing,
given in (\ref{reIII.1-las}).\\
(iv)\ For the subcase $\l=-1$, only $\varphi(L_n,Y_{-n-\mu})$ and
$\varphi(Y_n,M_{-n-3\mu})\,\,(\,\forall\,\,n\in\Z))$ are not
vanishing, which are given in (\ref{reIIILY1-1}) and
(\ref{corIIIYMa1-2}) respectively.
\end{lemm}
{\it Proof.} One has
\begin{eqnarray}
&&\varphi(L_0,Y_n)=0\ \ {\rm if}\ \,n\neq-\mu,\ \
\varphi(L_1,Y_{-\mu-1})=0\ \ {\rm if}\ \,\l\neq-3,\label{recIII2}\\
&&\varphi(L_0,M_n)=0\ \ {\rm if}\ \,n\neq-2\mu,\ \
\varphi(L_1,M_{-2\mu-1})=0\ \ {\rm if}\ \,\l\neq-1.\label{recIII3}
\end{eqnarray}
Using the Jacobian identity by replacing the triple $(v_1,v_2,v_3)$
by five triples $(Y_{m},\,M_n,\,L_0)$, $(M_{m},\,M_n,\,L_0)$,
$(Y_{m},\,Y_n,\,L_0)$, $(L_{m},\,Y_n,\,L_0)$ and
$(L_{m},\,M_n,\,L_0)$ in (\ref{2c-Ji}) respectively, one has
\begin{eqnarray}
&&(m+n+3\mu)\varphi(Y_m,M_n)=0,\label{corIIIYM}\\
&&(m+n+4\mu)\varphi(M_m,M_n)=0,\label{corIIIMM}\\
&&(m-n)\varphi(L_0,M_{m+n})+(m+n+2\mu)\varphi(Y_m,Y_n)=0,\label{corIIIYY}\\
&&(-n+\frac{\l+1}{2}m-\mu)\varphi(L_0,Y_{m+n})
+(m+n+\mu)\varphi(L_m,Y_n)=0,\label{corIIILY}\\
&&(-n+m\l-2\mu)\varphi(L_0,M_{m+n})+(m+n+2\mu)\varphi(L_m,M_n)=0.\label{corIIILM}
\end{eqnarray}
The following results can be directly obtained from
(\ref{corIIIYY})--(\ref{corIIILM}):
\begin{eqnarray}
&&\varphi(L_0,M_{-2\mu})=0\ \ \ {\rm and}\ \
(m+n+2\mu)\varphi(Y_m,Y_n)=0,\label{corIIIYYa1}\\
&&\varphi(L_0,M_{-2\mu})=0\ \ \ {\rm and}\ \
(m+n+2\mu)\varphi(L_m,M_n)=0,\label{corIIILMa1}\\
&&\varphi(L_0,Y_{-\mu})=0\ \ \ {\rm if}\ \ \l\ne-3,\ \ \ {\rm and}\
\ (m+n+\mu)\varphi(L_m,Y_n)=0.\label{corIIILYa1}
\end{eqnarray}
From (\ref{recIII2})--(\ref{recIII3}) and
(\ref{corIIIYYa1})--(\ref{corIIILYa1}), the left components we have
to present in this case are listed in the following (\,where $m$ is
an arbitrary integer):
\begin{eqnarray*}
&&\varphi(L_{-\mu-m},Y_m),\,\,\varphi(Y_{-2\mu-m},Y_m),\,\,
\varphi(L_{-2\mu-m},M_m),\,\,\varphi(Y_{-3\mu-m},M_m)\,\,{\rm
and}\,\,\varphi(M_{-2\mu-m},M_m),
\end{eqnarray*}
which will be taken into account in the following step by step.

\ni{\bf Step 1.} {\it The computation of
$\varphi(L_{-\mu-m},Y_m)\,\,(\,\forall\,\,m\in\Z)$}.

Replacing the triple $(v_1,v_2,v_3)$ by
$(L_{-\mu-m},\,L_{n},\,Y_{m-n})\,\,(\,\forall\,\,n\in\Z)$ in
(\ref{2c-Ji}), one has
\begin{eqnarray*}
&&\big(2(m+\mu)-n(\l+3)\big)\varphi(L_{-m-\mu},Y_m)\nonumber\\
&&=2(m+\mu+n)\varphi(L_{n-m-\mu},Y_{m-n})
+\big((m+\mu)(\l+3)-2n\big)\varphi(L_n,Y_{-n-\mu}),
\end{eqnarray*}
which gives (\,by replacing $n$ by $-n$ and $m$ by $m+n$
respectively)
\begin{eqnarray}
&&\big(2(m+\mu)+n(\l+3)\big)\varphi(L_{-m-\mu},Y_m)\nonumber\\
&&=2(m+\mu-n)\varphi(L_{-m-\mu-n},Y_{m+n})
+\big((m+\mu)(\l+3)+2n\big)\varphi(L_{-n},Y_{n-\mu}),\label{corIII.1ad-n}\\
&&\big(2(m+\mu+n)-n(\l+3)\big)\varphi(L_{-m-\mu-n},Y_{m+n})\nonumber\\
&&=2(m+\mu+2n)\varphi(L_{-m-\mu},Y_{m})
+\big((m+\mu+n)(\l+3)-2n\big)\varphi(L_n,Y_{-n-\mu}).\label{corIII.1adm+n}
\end{eqnarray}
If $\l\ne-3$, then taking $n=-1$ in both (\ref{corIII.1ad-n}) and
(\ref{corIII.1adm+n}), together with (\ref{recIII2}), one has
\begin{eqnarray}
&&\big(2(m+\mu)-(\l+3)\big)\varphi(L_{-m-\mu},Y_m)
=2(m+\mu+1)\varphi(L_{1-m-\mu},Y_{m-1}),\label{corIII.1ad-1}\\
&&\big(2(m+\mu-1)+(\l+3)\big)\varphi(L_{1-m-\mu},Y_{m-1})\nonumber\\
&&=2(m+\mu-2)\varphi(L_{-m-\mu},Y_{m})
+\big((m+\mu-1)(\l+3)+2\big)\varphi(L_{-1},Y_{1-\mu}).\label{corIII.1adm-1}
\end{eqnarray}
If $5-4\l-{\l}^2\ne0$, i.\,e., $\l\ne-5,1$, then combining
(\ref{corIII.1ad-1}) with (\ref{corIII.1adm-1}), we can deduce
\begin{eqnarray}\label{reIIILY1}
\varphi(L_{-m-\mu},Y_m)
=\frac{2(m+\mu+1)\big(2+(m+\mu-1)(3+\l)\big)\varphi(L_{-1},Y_{1-\mu})}
{5-4\l-{\l}^2}.
\end{eqnarray}
Taking $\varphi(L_{-m-\mu},Y_m)$ and $\varphi(L_{1-m-\mu},Y_{m-1})$
obtained from (\ref{reIIILY1}) back to (\ref{corIII.1ad-1}), one has
\begin{eqnarray*}
\frac{2\big((m+\mu)^2-1\big)(\l+1)(\l+3)
\varphi(L_{-1},Y_{1-\mu})}{5-4\l-{\l}^2}=0,
\end{eqnarray*}
which forces (since the index $m$ can be shifted and our assumption
$\l\ne-1$)
\begin{eqnarray*}
(\l+1)(\l+3) \varphi(L_{-1},Y_{1-\mu})=0.
\end{eqnarray*}
In another word, the system consisted of linear equations
(\ref{corIII.1ad-1}) and (\ref{corIII.1adm-1}) has nonzero solutions
if and only if $\l=-1$ or $\l=-3$ under our assumption
$5-4\l-{\l}^2\ne0$.

\ni If $\l=-1$, we can write (\ref{reIIILY1}) as
\begin{eqnarray}\label{reIIILY1-1}
\varphi(L_{-m-\mu},Y_m)
=\frac{(m+\mu)(m+\mu+1)}{2}\varphi(L_{-1},Y_{1-\mu}).
\end{eqnarray}
If $\l=-3$, we can write (\ref{reIIILY1}) as
\begin{eqnarray}\label{reIIILY1-2}
\varphi(L_{-m-\mu},Y_m) =\frac{m+\mu+1}{2}\varphi(L_{-1},Y_{1-\mu}).
\end{eqnarray}

If $\l=1$, then (\ref{corIII.1ad-1}) becomes
\begin{eqnarray}
&&(m+\mu-2)\varphi(L_{-m-\mu},Y_m)
=(m+\mu+1)\varphi(L_{1-m-\mu},Y_{m-1}),\label{corIII.1ad-1l=1}\\
&&(m+\mu+1)\varphi(L_{1-m-\mu},Y_{m-1})\nonumber\\
&&=(m+\mu-2)\varphi(L_{-m-\mu},Y_{m})
+(2m+2\mu-1)\varphi(L_{-1},Y_{1-\mu}),\label{corIII.1adm-1l=1}
\end{eqnarray}
which further gives \big(\,by taking $m=-\mu$ in
(\ref{corIII.1adm-1l=1})\big)
\begin{eqnarray}\label{corIIILY-ad-1}
\varphi(L_{-1},Y_{1-\mu})=0.
\end{eqnarray}
Also by (\ref{recIII2}), (\ref{corIIILYa1}), (\ref{corIII.1ad-1l=1})
and (\ref{corIIILY-ad-1}), one can deduce
\begin{eqnarray}\label{reIII.1-00a}
\varphi(L_{-m-\mu},Y_{m})=\left\{\begin{array}{ll}
(m+\mu-1)(m+\mu)(m+\mu+1)c_3&{\rm if}\ \,m\geq-\mu-1,\vs{6pt}\\
(m+\mu-1)(m+\mu)(m+\mu+1)c_4 &{\rm if}\ \,m<-\mu-1.
\end{array}\right.
\end{eqnarray}
where $c_3,\,c_4$ are some constants in $\C$. One thing left to be
done is to find the relations between the constants $c_3$ and $c_4$.
If $\l=1$, then (\ref{corIII.1ad-n}) and (\ref{corIII.1adm+n})
become
\begin{eqnarray*}
&&(m+\mu+2n)\varphi(L_{-m-\mu},Y_m)\nonumber\\
&&=(m+\mu-n)\varphi(L_{-m-\mu-n},Y_{m+n})
+(2m+2\mu+n)\varphi(L_{-n},Y_{n-\mu}),\\
&&(m+\mu-n)\varphi(L_{-m-\mu-n},Y_{m+n})\nonumber\\
&&=(m+\mu+2n)\varphi(L_{-m-\mu},Y_{m})
+(2m+2\mu+n)\varphi(L_n,Y_{-n-\mu}).
\end{eqnarray*}
which together force
\begin{eqnarray}
(2m+2\mu+n)
\big(\varphi(L_{-n},Y_{n-\mu})+\varphi(L_{n},Y_{-n-\mu})\big)=0,
\end{eqnarray}
and in particular give (by taking $n=2$)
\begin{eqnarray}
(m+\mu+1)
\big(\varphi(L_{-2},Y_{2-\mu})+\varphi(L_{2},Y_{-2-\mu})\big)=0.\label{reIII.1-1a}
\end{eqnarray}
Noticing $2-\mu>-1-\mu,\,-2-\mu<-1-\mu$, and then combining
(\ref{reIII.1-00a}) with (\ref{reIII.1-1a}), one can safely deduce
\begin{eqnarray}
c_3=c_4.\label{corIII.1adm+n2}
\end{eqnarray}
Hence (\ref{reIII.1-00a}) and (\ref{corIII.1adm+n2}) together give
\begin{eqnarray}\label{reIII.1-las}
\varphi(L_{-m-\mu},Y_{m}) =(m+\mu-1)(m+\mu)(m+\mu+1)c_3,\ \
\forall\,\,m\in\Z.
\end{eqnarray}

If $\l=-5$, then (\ref{corII.1adn}) and (\ref{corII.1adm+n}) convert
to the following form:
\begin{eqnarray}
&&(m+\mu-n)\varphi(L_{-m-\mu},Y_m)\nonumber\\
&&=(m+\mu-n)\varphi(L_{-m-\mu-n},Y_{m+n})
+(n-m-\mu)\varphi(L_{-n},Y_{n-\mu}),\label{corIII.1ad-5}\\
&&(m+\mu+2n)\varphi(L_{-m-\mu-n},Y_{m+n})\nonumber\\
&&=(m+\mu+2n)\varphi(L_{-m-\mu},Y_{m})
-(m+\mu+2n)\varphi(L_n,Y_{-n-\mu}).\label{corIII.1adm+5}
\end{eqnarray}
Furthermore, taking $n=-1$ in both (\ref{corIII.1ad-5}) and
(\ref{corIII.1adm+5}), and using (\ref{recIII2}), one has
\begin{eqnarray*}
(m+\mu+1)\big(\varphi(L_{-m-\mu},Y_{m})-
\varphi(L_{1-m-\mu},Y_{m-1})\big)\!\!\!&=\!\!\!&0,\\
(m+\mu-2)\big(\varphi(L_{-m-\mu},Y_{m})
-\varphi(L_{1-m-\mu},Y_{m-1})\big)\!\!\!&=\!\!\!&
(m+\mu-1)\varphi(L_{-1},Y_{1-\mu}),
\end{eqnarray*}
from which and using (\ref{recIII2}) again, can we deduce the
following relation:
\begin{eqnarray}\label{relalastIII0}
\varphi(L_{-m-\mu},Y_{m})=0,\ \ \forall\,\,m\in\Z.
\end{eqnarray}

\ni{\bf Step 2.} {\it The computations of $\varphi(L_{-2\mu-m},M_m)$
and $\varphi(Y_{-2\mu-m},Y_m),\,\,\forall\,\,m\in\Z$}.

Using the similar techniques used in Case \ref{cs2}, one also can
obtain the results listed below.

\ni If $\l=-2$, then
\begin{eqnarray}\label{relalast0V}
\varphi(L_{-m-2\mu},M_{m})=\varphi(Y_{-m-2\mu},Y_{m})=0,\ \
\forall\,\,m\in\Z.
\end{eqnarray}
If $\l=-1$, then
\begin{eqnarray}
&&\!\!\!\!\!\!\!\!\!\!
\varphi(L_{-m-2\mu},M_{m})=\left\{\begin{array}{ll}
(m+2\mu)\psi(L_1,M_{-1-2\mu}) &{\rm if}\ \,m\geq-2\mu,\vs{8pt}\\
(m+2\mu+2)\varphi(L_1,M_{-1-2\mu})\\
-(m+2\mu+1)\varphi(L_2,M_{-2-2\mu})&{\rm if}\ \,m<-2\mu,
\end{array}\right.\label{corIV.1ad-2-4}\\
&&\!\!\!\!\!\!\!\!\!\! \varphi(L_{-m-2\mu},Y_{m})
=\left\{\begin{array}{ll}
\big(\frac{m}{\mu}+1\big)\varphi(Y_{-2\mu},Y_0)\\
-m(\frac{m}{\mu}+2)\psi(L_1,M_{-1-2\mu}) &{\rm if}\ \,m\geq0,\vs{8pt}\\
\big(\frac{m}{\mu}+1\big) \varphi(Y_{-2\mu},Y_0)
-(\frac{m}{\mu}+2)\big((m+2)\\
\times\varphi(L_1,M_{-1-2\mu})-(m+1)\varphi(L_2,M_{-2-2\mu})\big)&{\rm
if}\ \,m<0.
\end{array}\right.\label{N.2ad-2}
\end{eqnarray}
\ni If $\l=0$, then
\begin{eqnarray}
&&\!\!\!\!\!\!\varphi(L_{-m-2\mu},M_m)\nonumber\\
&&\!\!\!\!\!\!=\frac{(m+2\mu)(m+2\mu+1)}{2}
\varphi(L_{-1},M_{1-2\mu}),\ \ \ \forall\,\,m\in\Z;\label{reIV1l=066}\\
&&\!\!\!\!\!\!\varphi(Y_{-4\mu},Y_{2\mu})=\frac{5}{3}\varphi(Y_{-2\mu},Y_{0})
+\frac{4\mu(4\mu+1)}{3}\varphi(L_{-1},M_{1-2\mu}),\label{N.2nesp2-166}\\
&&\!\!\!\!\!\!\varphi(Y_{-m-2\mu},Y_{m})\nonumber\\
&&\!\!\!\!\!\!=\frac{m(m+1)(m+2\mu)}{m-2\mu}
\varphi(L_{-1},M_{1-2\mu})-\frac{3m+2\mu}{m-2\mu}\varphi(Y_{-2\mu},Y_0),\
\ {\rm if}\,\,m\ne2\mu.\label{N.2lad-3}
\end{eqnarray}
\ni If $\l=1$, then
\begin{eqnarray}
&&\varphi(Y_{-\mu-2\mu},Y_{\mu})
=\frac{1}{4}\mu(\mu^2-1)c'_1-\frac{5}{4}\varphi(Y_{0},Y_{-2\mu}),\label{N.2lm=mu}\\
&&\varphi(Y_{-m-2\mu},Y_{m})=\frac{2m+\mu}{\mu-m}\varphi(Y_{-2\mu},Y_0),\
\ \ {\rm if}\,\,m\ne\mu;\label{N.2lad-5}\\
&&\varphi(L_{-m-2\mu},M_{m}) =(m+2\mu-1)(m+2\mu)(m+2\mu+1)c'_1,\ \
\forall\,\,m\in\Z.\label{reIV.1-las}
\end{eqnarray}
where $c'_1$ is some constant in $\C$. If
$\l\notin\{-2,\,-1,\,0,\,1\}$, then
$\varphi(L_{-m},M_{m-2\mu})=0,\,\,\forall\,\,m\in\Z$.

\ni If $\frac{2\mu}{\l+1}\notin\Z$, then
\begin{eqnarray}\label{N.2lad-6-1}
&&\varphi(Y_{-m-2\mu},Y_{m})
=-\frac{m(\l+3)+2\mu}{m(\l+1)-2\mu}\varphi(Y_{-2\mu},Y_0).
\end{eqnarray}
If $\frac{2\mu}{\l+1}\in\Z,\,\l\notin\{-3,\,-2,\,-1,\,0,\,1\}$, then
\begin{eqnarray}\label{N.2lad-6-1}
&&\varphi(Y_{-m-2\mu},Y_{m})=0,\ \ \ \forall\,\,m\in\Z.
\end{eqnarray}
If $\l=-3$, then
\begin{eqnarray}\label{N.2nesp2las6}
\varphi(Y_{-m-2\mu},Y_{m})=\left\{\begin{array}{cl}
0 &{\rm if}\ \,m\geq\mu,\vs{8pt}\\
-\frac{1}{2}\varphi(Y_{0},Y_{-2\mu})&{\rm if}\ \,m=\mu.
\end{array}\right.
\end{eqnarray}

\ni{\bf Step 3.} {\it The computation of
$\varphi(Y_{-3\mu-m},M_m),\,\,\forall\,\,m\in\Z$}.

Replacing the triple $(v_1,v_2,v_3)$ by
$(L_{-2\mu-m},\,Y_{-\mu},\,M_m)$ in (\ref{2c-Ji}), one has
\begin{eqnarray*}
(m+2\mu)(1+\l)\big(\varphi(Y_{-m-3\mu},M_m)
+2\varphi(Y_{-\mu},M_{-2\mu})\big)=0,
\end{eqnarray*}
which gives, if $\l\ne-1$,
\begin{eqnarray}
\varphi(Y_{-m-3\mu},M_m)=\left\{\begin{array}{ll}
\varphi(Y_{-\mu},M_{-2\mu})&{\rm if}\ \,\, m=-2\mu,\vs{6pt} \\
-2\varphi(Y_{-\mu},M_{-2\mu})&{\rm if}\ \,\, m\ne-2\mu.
\end{array}\right.
\end{eqnarray}
Replacing the triple $(v_1,v_2,v_3)$ by $(L_{-m},\,Y_{-3\mu},\,M_m)$
in (\ref{2c-Ji}), one has
\begin{eqnarray*}
\big(m(1+\l)-4\mu\big)\varphi(Y_{-m-3\mu},M_m)=-2\big(m(1+\l)+2\mu\big)
\varphi(Y_{-3\mu},M_0),
\end{eqnarray*}
which can be rewritten as follows, if $\l=-1$,
\begin{eqnarray}
\varphi(Y_{-m-3\mu},M_m)=\varphi(Y_{-3\mu},M_0)\ \ \ ({\rm as}\
\mu\in\Z^*\ {\mbox{in this case}}).\label{corIIIYMa1-2}
\end{eqnarray}

\ni{\bf Step 4.} {\it The computation of
$\varphi(M_{-4\mu-m},M_m),\,\,\forall\,\,m\in\Z$}.

Finally, similar as that of Subcase 1.1. we can prove
\begin{eqnarray}
&&\varphi(M_{-4\mu-m},M_m)=0,\ \
\forall\,\,m\in\Z.\label{corIII-I-1}
\end{eqnarray}

Since for any $x,y,z\in\LL$, the Jacobian identity must be
satisfied, we can obtain all the compatible cocycles. They are just
those listed in Lemma \ref{lemm4} in this case. This completes the
proof of the lemma.\QED

\noindent Then Theorem \ref{mainth} (iv)--(vi) follow from Lemma
\ref{lemm4}. And above all, the main theorem can be easily deduced
from lemma \ref{lemm1}--\ref{lemm4}.


\begin{thebibliography}{9999}\vskip0pt\small
\parindent=2ex\parskip=-1pt\baselineskip=-1pt

\bibitem{BKV} B.~Bakalov, V.G.~Kac, A.A.~Voronov, ``Cohomology of
conformal algebras,'' {\it Comm.~Math. Phys.} {\bf200} (1999),
561--598.\\[-10pt]
\bibitem{H1} M. Henkel, Schr\"{o}dinger invariance and strongly
anisotropic critical systems, {\it J. Stat. Phys.}, {\bf 75} (1994),
1023-1029.\\[-10pt]
\bibitem{H2} M. Henkel, Phenomenology of local scale invariance: from
conformal invariance to dynamical scaling, {\it Nucl. Phys. B}, {\bf
641} (2002), 405-410.\\[-10pt]
\bibitem{HU} M. Henkel, J. Unterberger, Schr\"{o}dinger invariance and
space-time symmetries, {\it Nucl. Phys. B}, {\bf 660} (2003),
407-412.\\[-10pt]
\bibitem{LS1} J. Li, Y. Su, The derivation algebra and automorphism
group of the twisted Schr\"{o}dinger-Virasoro algebra, preprint.\\[-10pt]
\bibitem{LS2} J. Li, Y. Su, Representations of the Schr\"{o}dinger-Virasoro
algebras, preprint.\\[-10pt]
\bibitem{L}   W.~Li, ``2-Cocycles on the algebra of differential
operators,'' {\it J.~Algebra} {\bf122} (1989), 64--80.\\[-10pt]
\bibitem{LW}   W.~Li, R.L.~Wilson, ``Central extensions of some Lie
algebras,'' {\it Proc.~Amer.~Math.~Soc.} {\bf126} (1998),
2569--2577.\\[-10pt]
\bibitem{RU} C. Roger, J. Unterberger, The Schr\"{o}dinger-Virasoro Lie
group and algebra: from geometry to representation theory, preprint
(arXiv:cond-mat/0601050), (2006).\\[-10pt]
\bibitem{SZR} M.~Scheunert, R.B.~Zhang, ``Cohomology of Lie superalgebras and
their generalizations,'' {\it J.~Math. Phys.} {\bf39} (1998),
5024--5061.\\[-10pt]
\bibitem{S1} Y.~Su, ``2-Cocycles on the Lie algebras of all differential
operators of several indeterminates,'' {\it (Chinese) Northeastern
Math.~J.} {\bf6} (1990), 365--368.\\[-10pt]
\bibitem{S2} Y.~Su, ``2-cocycles on the Lie algebras of generalized
differential operators'', {\it Comm.~Algebra} {\bf30} (2002),
763--782.\\[-10pt]
\bibitem{S3} Y.~Su, ``Low dimensional cohomology of general conformal algebras
$gc_N$,'' {\it J.~Math.~Phys.} {\bf45} (2004), 509--524.\\[-10pt]
\bibitem{SZK}  Y.~Su, K.~Zhao, ``Second cohomology group of generalized Witt
type Lie algebras and certain reperesentations,'' {\it
Comm.~Alegrba} {\bf30} (2002),  3285--3309.\\[-10pt]
\bibitem{U} J. Unterberger, On vertex algebra representations of the
Schr\"{o}dinger-Virasoro algebra, preprint (arXiv:cond-mat/0703214),
(2007).








\end{thebibliography}
\end{document}